\documentclass[12pt]{article}

\usepackage{amsmath}
\DeclareMathOperator*{\argmin}{arg\,min}
\usepackage{amsfonts}
\usepackage{amsthm}
\newtheorem{claim}{Claim}
\usepackage{url}
\usepackage{times}
\usepackage{graphicx}
\usepackage{subcaption}
\usepackage[ruled,vlined]{algorithm2e}

\newenvironment{sciabstract}{%
\begin{quote} \bf}
{\end{quote}}
\newcommand*{\vertbar}{\rule[-1ex]{0.5pt}{2.5ex}}

\title{DySMHO: Data-Driven Discovery of Governing Equations for Dynamical Systems via Moving Horizon Optimization}

\author
{Fernando Lejarza,$^{1}$ Michael Baldea,$^{1,2\ast}$\\
\\
\normalsize{$^{1}$McKetta Department of Chemical Engineering, The University of Texas at Austin,}\\
\normalsize{$^{2}$ODEN Institute for Computational Engineering and Sciences,The University of Texas at Austin}\\
\\
\normalsize{$^\ast$To whom correspondence should be addressed; E-mail:  \url{mbaldea@utexas.che.edu}}
}

\date{}

\begin{document}

\newcommand{\ours}{DySMHO}

\baselineskip18pt

\maketitle

\section*{Abstract}
\begin{sciabstract}
Discovering the governing laws underpinning physical and chemical phenomena is a  key step towards understanding and ultimately controlling  systems in  science and engineering. We introduce Discovery of Dynamical Systems via Moving Horizon Optimization (DySMHO), a scalable machine learning framework for identifying governing laws in the form of differential equations from large-scale noisy experimental data sets. \ours{} consists of a novel moving horizon dynamic optimization strategy that sequentially  learns the underlying governing equations from a large dictionary of basis functions. The sequential nature of DySMHO allows leveraging statistical arguments for eliminating irrelevant basis functions, avoiding overfitting to recover accurate and parsimonious forms of the governing equations. Canonical nonlinear dynamical system examples are used to demonstrate that \ours{} can accurately recover the governing laws, is robust to high levels of measurement noise and that it can handle challenges such as multiple time scale dynamics.
\end{sciabstract}

\section*{Introduction}
Differential and partial differential equation models play a critical role in describing the governing behavior of a variety of systems arising in science and engineering \cite{gockenbach2005partial}. As minimal order expressions describing the system behavior, governing models are generalizable and readily interpretable, and have  good extrapolation capabilities. Historically, the discovery and formulation of fundamental governing equations is a relatively lengthy process,  supported by careful experimentation and data collection using prototype systems.

The data sets that have, through the history of science, supported the discovery of fundamental natural laws may seem ``small'' by today's standards. With decreasing costs of sensors, data storage systems, and computing hardware, immense quantities of data can be easily collected and efficiently stored. As a consequence, the applications of machine learning (ML) and artificial intelligence (AI) have witnessed meteoric growth in science and engineering. ML techniques perform well in regression and classification tasks, but the resulting models are   ``black'' or ``grey-box'' in nature, offering little physical insight \cite{roscher2020explainable}, and extrapolate poorly to regimes beyond the scope of the training data. Moreover, prediction accuracy typically comes at the cost of model complexity, which is at odds with the parsimonious  nature of a system's governing dynamics derived via first principles analysis.

Leveraging ML/AI frameworks to discover (as opposed to merely fit) the governing equations of physical systems from large amounts of data offers intriguing possibilities and remains an open field of research. Recent efforts in this direction include physics-informed discovery strategies \cite{karniadakis2021physics}, combining first principles arguments with ML models such as Gaussian processes \cite{raissi2018hidden} and deep neural networks \cite{raissi2019physics,long2019pde}. While the predictive capabilities of such physics-informed strategies are good, even when trained on coarse and noisy measurements, their reliance on significant (if not complete) structural knowledge of the equations governing the system dynamics constitutes a significant disadvantage when it comes to discovering they governing equations of new and unknown systems. Further, as is inherent to most (deep) machine learning architectures, such models can suffer from lack of interpretability, hence failing to provide insights on the selection process of the functional terms that dictate a system's dynamic behavior.

A  different approach that is generally deemed more transparent  towards automating the data-driven discovery of governing equations is based on nonlinear regression strategies \cite{james2013introduction}. Initial efforts exploited symbolic regression  \cite{king2004functional,bongard2007automated, schmidt2009distilling} and genetic programming algorithms \cite{koza1992genetic}. However, the combinatorial nature of these approaches can render them computationally prohibitive, restricting their applicability to low-dimensional systems and to considering relatively small initial sets of candidate symbolic expressions (which inherently diminishes the probability of identifying the true underlying system dynamics). Furthermore, symbolic regression strategies are prone to overfitting, i.e., generating overly complex expressions in an attempt to decrease prediction error \cite{schmidt2009distilling}.

More recently, sparse regression techniques \cite{tibshirani1996regression,zou2005regularization} have been proposed. Brunton et al. \cite{brunton2016discovering} employed a modified ordinary least-squares (OLS) and LASSO regression (i.e., $\ell_1$ penalized regression) to discover parsimonious representations of the dynamics of nonlinear systems  from high-dimensional data sets by selecting the elements of the governing equations form a \textit{a priori} specified large set of candidate basis functions. Numerous extensions to \cite{brunton2016discovering} have since been proposed addressing a variety of classes of systems and problem settings (e.g., \cite{mangan2016inferring,mangan2017model,kaiser2018sparse,champion2019data}).
A similar approach based on elastic net regression (i.e., a combination of both $\ell_1$ and $\ell_2$  penalized regression) was introduced in \cite{sun2020alven}, but resulted in less parsimonious equations relative to e.g. LASSO regression. In a related effort \cite{cozad2014learning}, the selection of basis functions was performed via mixed-integer optimization, with a view towards identifying low-order surrogate representations of nonlinear algebraic models. 

In spite of these advances, several  fundamental challenges  remain. From a numerical perspective, existing approaches rely on either directly measuring the system state derivatives (which are likely not observable in practical settings) or approximating them accurately (which can be  difficult particularly in high noise environments and when the dynamics evolve over multiple time scales). Additionally, when considering large libraries of nonlinear basis functions, feature collinearity may result in ill-conditioned regression problems. As the size of the available data set increases, the ensuing numerical instability renders these approaches unrealiable for discovering fundamental equations that are optimal (in the Pareto sense) with respect to both model parsimony and predictive power.

Motivated by the above, in this work, we propose \ours, a radically different perspective to discovering governing equations from data. The present method rooted in control theory, namely, moving horizon estimation and control, and offers (i) excellent scalability with respect to the system dimensions and the size of the data set, (ii) rigorous statistical arguments for selecting the model structure from a large dictionary of basis functions, (iii) robustness to noisy training data and (iv) the ability to incorporate first principles knowledge (when available) in the form of additional constraints in the problem formulation.

\section*{Results}

\subsection*{Representation of system dynamics}
We consider dynamical systems governed by ordinary differential equations of the form:
\begin{equation}
\label{eq:dyn_system}
\frac{d}{dt}\textbf{x}(t) = \textbf{f}(\textbf{x}(t))
\end{equation}
where $\textbf{x}(t) \in \mathbb{R}^{n_x}$ is the vector of states at time $t$, and the map $\textbf{f}(\cdot):\mathbb{R}^{n_x}\rightarrow \mathbb{R}^{n_x}$ represents the (nonlinear) dynamics of  the system. The function $\textbf{f}$ is \textit{unknown} and is precisely what we attempt  to infer from a given set of time-resolved measurement data. To that end, we collect a sequence of measurements $\hat{\textbf{x}}(t_k)$ of the state variables observed at sampling times $t_1,\dots,t_m$, and assume that the derivative $\dot{\hat{\textbf{x}}}(t_k)$ cannot be directly observed. The data are assumed to be contaminated with (Gaussian, zero-mean) measurement noise, and smoothing techniques and statistical tests are used to perform pre-processing  of the training data set (See Materials and Methods). The resulting pre-processed data are denoted by $\tilde{\textbf{x}}(t_k) \; \forall k=1,\dots,m$.

To discover the underlying governing equations, we consider a dictionary of $n_\theta$ candidate symbolic nonlinear basis functions denoted as $\Theta(\textbf{x}^T)$, where $ \Theta(\cdot): \mathbb{R}^{1\times n_x} \rightarrow \mathbb{R}^{1\times n_\theta}$. The dictionary is defined \textit{a priori}, potentially leveraging some domain insights (e.g. \cite{sun2020alven}). Importantly, we assume that the governing equations can be expressed as an (as of now unknown) linear combination of the basis functions in this dictionary or, equivalently, that the true model is contained within the dictionary. That is, we can express  the model in \eqref{eq:dyn_system} as:
\begin{equation}
\label{eq:dyn_basis}
\frac{d}{dt}\textbf{x}(t) = \textbf{f}(\textbf{x}(t)) = \Xi^T (\Theta (\textbf{x}(t)^T))^T
\end{equation}
where $\Xi \in \mathbb{R}^{n_\theta \times n_x}$ is a matrix of coefficients whose columns are given by the sparse (i.e., most entries are zero) vectors $\pmb{\xi}_1, \dots,\pmb{\xi}_{n_x}$. This sparse model structure is illustrated in Fig. 1 (A) using the well-known two-dimensional Lotka-Volterra predator-prey model as an example. Further, we distinguish between \textit{basic} coefficients (i.e., $\xi_{i,j} \neq 0$ in the true dynamics, and the corresponding basis functions are active), and \textit{non-basic} coefficients (i.e., $\xi_{i,j}=0$ in the true dynamics). Basic coefficients within $\Xi^T$ for the predator-prey model are highlighted in color in Fig. 1 (A). Following the sparsity argument, the fundamental challenge of discovering the governing equations translates to identifying the (few) basic coefficients associated with the active nonlinearities that compose $\textbf{f}(\textbf{x}(t))$. Fig. 1 (B) shows first an example dictionary of basis functions evaluated over time for the Lotka-Volterra model, where the active basis functions relevant to describing the dynamics of  each of the two system states are shown in color. A linear combination of the selected basis functions is then used to construct the  dynamics for each state.

\begin{figure}[h!]
\centering
\includegraphics[scale = 0.3]{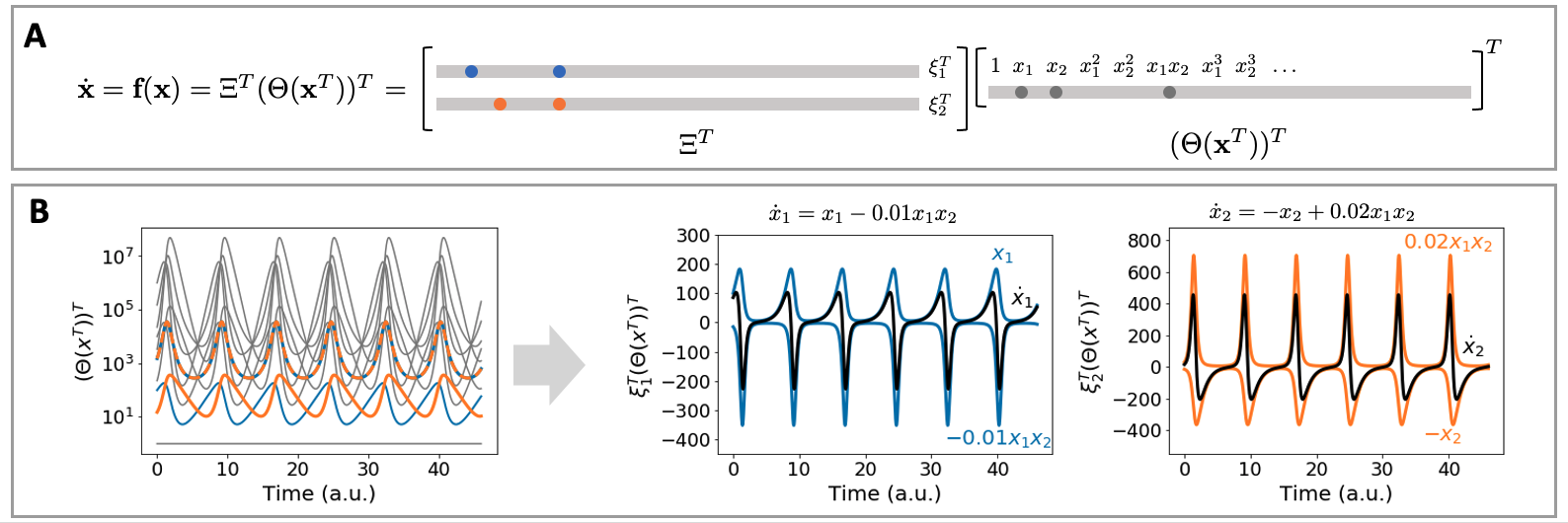}
\caption{\ours{} model structure for discovering the the Lotka-Volterra predator-prey model, with $\dot x_1 = x_1 - 0.01 x_1 x_2$, $\dot x_2 = -x_2 + 0.02 x_1 x_2$. (\textbf{A}) Graphical representation of the dictionary of basis functions $\Theta(\textbf{x}^T)$  and sparse coefficient matrix $\Xi$ in two dimensions. (\textbf{B}) Illustration of the values of basis functions for each of the two states in predator-prey model (colored lines indicate the true basis functions against the backdrop of all basis functions which are shown in grey).}
\label{fig:model_structure}
\end{figure}

Prior works (e.g., \cite{brunton2016discovering,sun2020alven}) proposed sparse regression techniques, whereby  for each  state variable $i \in \{1,\dots,n\}$ the sparse vector $\pmb{\xi}_i$ is estimated as the solution of the following optimization problem:
\begin{equation}
\label{eq:regularized_regression}
\pmb{\xi}_i \in   \argmin \frac{1}{2m} \sum_{j=1}^{m}||\Theta(\tilde{\textbf{x}}(t_j)^T)\pmb{\xi}_i - \dot{\tilde{\textbf{x}}}(t_j) ||_2^2 + \lambda \rho ||\pmb{\xi}_i||_1 + \frac{\lambda(1-\rho)}{2} ||\pmb{\xi}_i||_2^2
\end{equation}
The goal is to minimize, in a norm sense, the difference between the value of the state derivatives predicted by the model, $\Theta(\tilde{\textbf{x}}(t_j)^T)\pmb{\xi}_i$, and the corresponding  derivative values $\dot{\tilde{\textbf{x}}}(t_j)$ that are either measured  directly or approximated  (via e.g. finite difference equations)  from the $m$ data samples at each sample time $t_j$. The second part of the expression above is a regularization term that places a penalty on the   $\ell_1$ and/or  $\ell_2$-norms of the magnitudes of the coefficients $\pmb{\xi}_i$, thereby ensuring that the solution is sparse (i.e., that as many of the elements of this vector are zero).  $\lambda \ge 0$ controls the extent of this spasification, while  $0 \le \rho < 1$ allows for balancing between $\ell_1$ and/or  $\ell_2$ norms.

\subsection*{Nonlinear optimization-based discovery formulation}

Differently, to obtain the sparse coefficient matrix $\Xi$ in \eqref{eq:dyn_basis}, we formulate a constrained multi-period dynamic nonlinear program (DNLP) minimizing a given error metric. Discretization strategies are used to convert the candidate model from the differential form in \eqref{eq:dyn_basis} to a system of nonlinear algebraic equations, which is embedded in  the constraint system of  the DNLP \cite{nicholson2018pyomo} (See Materials and Methods for further details). Pre-processing techniques and/or domain expertise can be used to introduce additional constraints to the DNLP in order to improve convergence to the governing equations (See Materials and Methods for further details). The DNLP is thus formulated in discrete time to minimize the mean $\ell_2$-error as follows:
\begin{equation}
\label{eq:dyn_opt}
\begin{aligned}
& \min_{\Xi} && \frac{1}{2M} \sum_{k=1}^{M}||\textbf{x}(k)-\tilde{\textbf{x}}(k)||_2^2 + \lambda \ell(\Xi)\\
& \; \text{s.t.} && \textbf{x}({k+1}) = \textbf{g}(\Theta(\textbf{x}(k)),\Xi)  \;\;\; \forall k\in \{1,\dots,M\}\\
& \; &&\textbf{x}(1) = \tilde{\textbf{x}}(t_1) \\
& \; && \Xi \in \{ \Xi^{L}, \Xi^{U} \}
\end{aligned}
\end{equation}
where $\textbf{x}(k)$ are the states predicted by the candidate governing law at time index $k$, $\textbf{g}(\cdot):\mathbb{R}^{n_x}\times \mathbb{R}^{n_\theta \times n_x}\rightarrow\mathbb{R}^{n_x}$ represents the discretized version of the nonlinear dynamics given in \eqref{eq:dyn_basis}, $\Xi^{L}$ and $\Xi^{U}$ are respectively the lower and upper bounds of the estimated coefficients, and $M$ is the number of discrete time points of the transformed dynamics. Note that in \eqref{eq:dyn_opt} the basis functions in the dictionary are not directly evaluated on the measurement data $\tilde{\textbf{x}}$ as is the case in \eqref{eq:regularized_regression}. Rather,  expression \eqref{eq:dyn_opt} represents a symbolic nonlinear function of the (predicted) states, which are decision variables in the DNLP. Further, \eqref{eq:dyn_opt} does not (directly) require an approximation of $\tilde{\dot{\textbf{x}}}$ to derive the coefficients $\Xi$  (as is done in the sparse regression problem \eqref{eq:regularized_regression}). The derivative approximation $\dot{\tilde{\textbf{x}}}$ is used indirectly in pre-processing (as discussed in Materials and Methods) as a potential first step for pruning the basis function dictionary, as well as for estimating $\Xi^L$ and $\Xi^U$.

An additional advantage of the proposed formulation is that the dynamics of  \textit{all} states are simultaneously recovered, as opposed to prior works \cite{brunton2016discovering, sun2020alven} in which regression problems are solved \textit{separately} for each measured state. It is thus expected that the solution of \eqref{eq:dyn_opt} results in governing equation models that better explain the overall system behavior, avoiding possible additional variance that might be introduced when the dynamics are independently identified for each observed state variable. Note that interpolation is typically required to approximate the state measurements $\tilde{\textbf{x}}(t_k)\; \forall k=1,\dots,m$ at each of the points in the new time grid defined by the discretization method employed to construct $\textbf{g}(\cdot)$. The objective function in \eqref{eq:dyn_opt} includes a regularization term, where the parameter $\lambda$ penalizes the norm of the coefficient vector $\Xi$ and imposes sparsity in the same sense as described earlier.

The dimension  of the DNLP \eqref{eq:dyn_opt} increases with the dimension $M$ of the data set available and with the dimension $n_\theta$  of the dictionary of basis functions. Intuitively, $M$ may be very large (i.e., many data points are available from experiments), while a large $n_\theta$ (i.e., a large dictionary of basis functions) is highly desirable in order to increase the probability of discovering the true governing equations \eqref{eq:dyn_system}. Thus, solving \eqref{eq:dyn_opt} while taking into consideration the entire available data set, as is generally done in existing discovery frameworks based on sparse regression\cite{brunton2016discovering,sun2020alven}, is likely computationally expensive. Noting that such problems are NP-hard, the actual solution time cannot be anticipated from the above problem dimensions. An additional fundamental challenge is related to imposing parsimony in the learned model. This entails eliminating the basis functions that are not part of the true model \eqref{eq:dyn_system} by setting the corresponding coefficients $\Xi$ to zero. A particular difficulty arises when the true value of a coefficient is ``small:'' while the corresponding estimate may also be small, it is difficult to discern whether this outcome is correct or the non-zero estimated value is the result of ovefitting (i.e., a spurious attempt to further decrease the value of the objective function in \eqref{eq:dyn_opt} by increasing model complexity). A thresholding approach consisting of eliminating terms whose estimated coefficients are below a specific value (determined via cross-validation) can in principle be employed, but its performance is expected to degrade with increasing model stiffness and to our knowledge there is no rigorous way of defining such thresholds other than through cross-validation.

The framework proposed here addresses both fundamental problems described above. To deal with problem dimensionality,  we propose decomposing \eqref{eq:dyn_opt} into a sequence of lower-dimensional problems defined on shorter time horizons (i.e., using smaller subsets of the available data), for which optimal solutions to \eqref{eq:dyn_opt} can be attained with significantly lower computational effort. An illustration of the data subsets is shown in Fig. S.2 in the Supporting Materials. After a solution to \eqref{eq:dyn_opt} is computed, a new data subset is selected (intuitively -- but not necessarily -- by shifting the time window forward by a smaller step than the window size) which is then used to resolve \eqref{eq:dyn_opt} again. The repetition of this procedure allows for efficiently learning and refining a sequence of governing equation models each with different coefficient estimates (further details regarding the moving horizon optimization and thresholding strategies proposed can be found in Materials and Methods and in the Supporting Materials). In conjunction with this moving-horizon strategy, the following thresholding claim is made: 

\begin{claim}
The non-basic coefficients in $\Xi$  typically contribute to overfitting. Small, non-zero values reflect the use of the corresponding functions to fit the noise in the training data. Hence, the mean of the value of the estimates of a non-basic coefficient derived from the sequence of problems described above is likely to have a relatively high standard deviation. 
\end{claim}

The converse  argument can be made for basic coefficients: the variance of a sequence of estimates is expected to be relatively low. These claims then support the use of dispersion metrics in statistics (e.g., the coefficient of variation) for the parameters obtained in a sequence of estimates based on subsets of the data to infer whether a parameter belongs to the true dynamics or not (i.e., if it is basic or non-basic). Note the similarity of the proposed training mechanism with ensemble methods (e.g., \cite{ho1995random}) where a pool of models is trained on different (typically randomly sampled) subsets of the data to reduce variance and thus improve the performance of the final model prediction. Our moving horizon and thresholding mechanisms are discussed in depth in the Materials and Methods sections, and detailed algorithm steps are provided in the Supporting Materials. A comprehensive illustration of the entire workflow of the \ours{} framework is shown in Fig. 2.

\begin{figure}[h!]
\centering
\includegraphics[scale = 0.3]{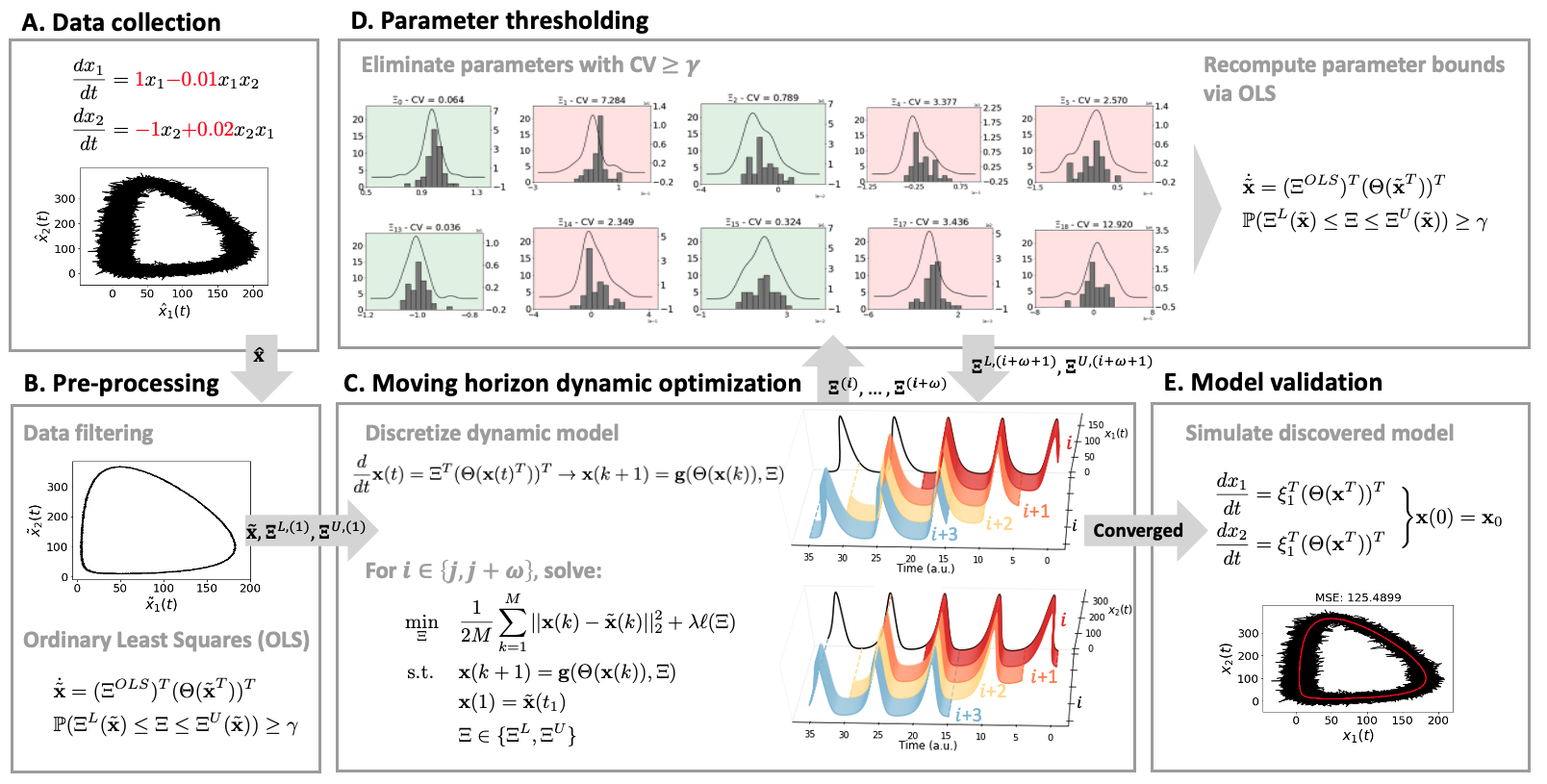}
\caption{Illustration of the \ours{} workflow: Noisy data $\hat{\textbf{x}}$ are initially collected. The data are smoothed and pre-processed and bounds on basis function coefficients are established. The smoothed data $\tilde{\textbf{x}}$ and the bounds  ($\Xi^{L,(i)}$, $\Xi^{U,(i)}$) are used to formulate the moving horizon DNLP. Every $\omega$ iterations, thresholding is performed to eliminate non-basic parameters and coefficient bounds are tightened. Convergence is achieved once the set of basic functions does not change for a given number of iterations, or once the entire training data set has been exhausted. The discovered model is validated via simulation and comparison with the data. Qualitative comparisons can be performed using, e.g.,phase plane plots, while quantitative assessments rely on common regression performance metrics (e.g., mean squared error).}
\label{fig:schematic_v1}
\end{figure}

To evaluate the performance of \ours,  we consider a series of canonical nonlinear dynamical systems: the Lotka-Volterra predator-prey model \cite{lotka1925elements,volterra1926fluctuations}, the van der Pol oscillator \cite{van1926lxxxviii}, the Brusselator \cite{prigogine1968symmetry}, and the chaotic Lorenz oscillator \cite{lorenz1963deterministic}. In each case, the true model  is used to generate data via simulation, and data are artificially contaminated with  Gaussian noise of increasing standard deviation $\sigma$. The results from our numerical experiments are shown in Fig. 3 and Fig. 4, where we explore the convergence of \ours{} in relation to the discovered model's complexity (i.e., number of terms), accuracy (i.e., error in the coefficient estimates), and predictive capabilities (i.e., mean squared error with respect to the true system dynamics). Complete details regarding the computational experiments performed as well as the numerical parameters and configurations for \ours{} are reported in the Supporting Materials Tables S.1-S.4.

\begin{figure}[h!]
\includegraphics[scale = 0.4]{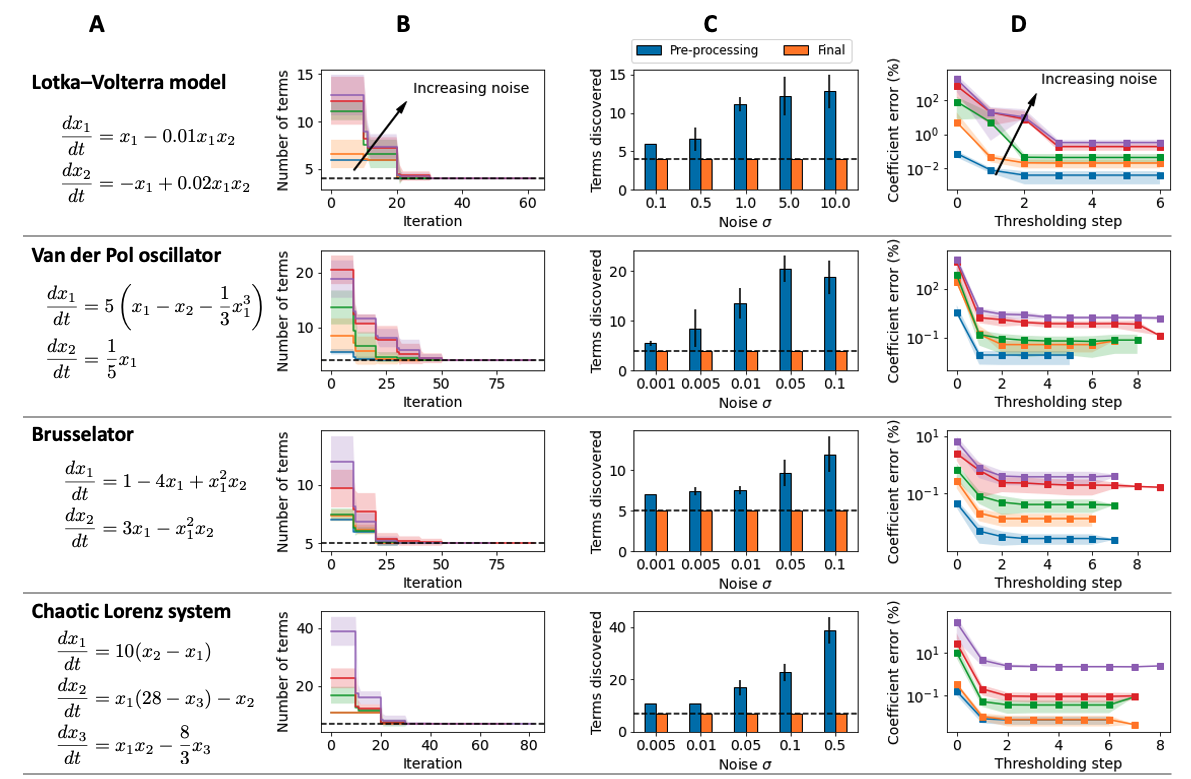}
\caption{Numerical experiments of \ours{} for the dynamical systems considered. (\textbf{A}) True governing differential equations for dynamical systems. (\textbf{B}) Average number of basis functions remaining in the discovered model as a function of \ours{} iterations; a dashed black line is used to indicate the true number of terms (two-dimensional systems were initially modeled with 28 basis functions, and 3-dimensional systems with 66 (See Supporting Materials for further details)). (\textbf{C}) Average number of terms in the discovered governing equation after pre-processing and the final model after \ours{} converged. (\textbf{D}) Percent cumulative coefficient error relative to the coefficients in the true governing equations after each thresholding step for data simulated with increasing noise contamination. (Results show mean values obtained from 10 random samples of the simulated measurement data for each noise level considered, and error bars and shaded areas represent one standard deviation from the mean)}
\end{figure}

\begin{figure}[h!]
\includegraphics[scale = 0.4]{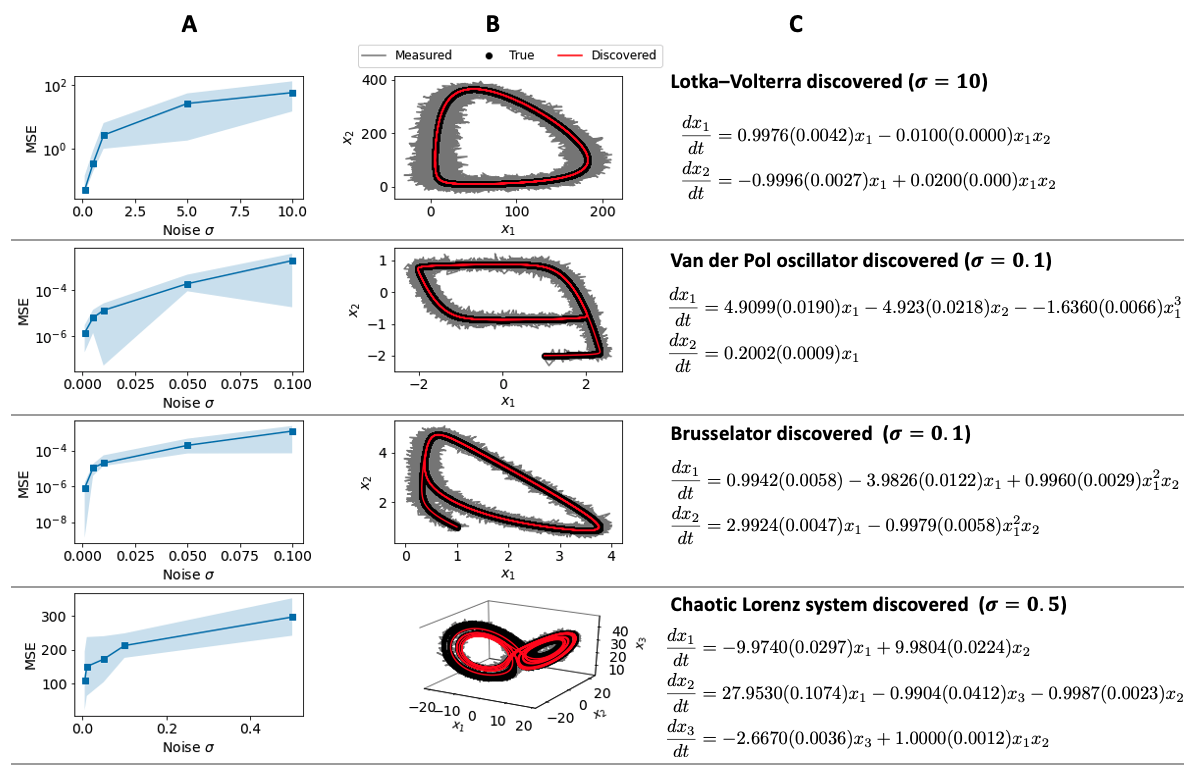}
\caption{Validation of the discovered equations obtained through \ours{} for the dynamical systems considered. (\textbf{A}) Average mean squared error (MSE) between simulation obtained from the true model and simulation obtained for the discovered model for 10 realizations of training data for each noise level. (\textbf{B}) Comparison of discovered dynamics against true model and measured data (data measured with $\sigma =10, 0.1, 0.1, 0.5$ respectively for Lotka-Volterra model, van der Pol oscillator, Brusselator, Lorenz system. (\textbf{C}) Discovered differential equations showing mean coefficient estimates (with standard deviation shown in parentheses) for the highest noise level considered for each dynamical system.}
\end{figure}

\section*{Discussion}

\ours{} accurately recovers the true governing equations for all systems, even in cases where data are noisy (for noise levels up to the standard deviations shown in Fig. 3 and Fig. 4). Discovery for higher noise settings than the ones considered can be accomplished by tuning pre-processing and moving horizon parameters on a system-by-system basis (e.g. the thresholding tolerance for the coefficient of variation of a sequence of model coefficient estimates).  Evidently, as shown in Fig. 4 (B) as the amount of measurement noise increases, so does the mean squared error (MSE) between the simulated state trajectories from the recovered equations and that of the simulation based on the true dynamics. The higher observed MSE values correspond to higher coefficient error (relative to the true coefficients) as shown in Fig. 3 (D). While the dynamics are seen to converge to the correct functional form of the true governing dynamics, noisier data inherently introduce additional uncertainty in the values of the coefficient estimates. Under significant coefficient estimate variability, Monte Carlo sampling can be used to generate candidate models from which several plausible state trajectories can be evaluated to quantify uncertainty in the discovered dynamics.

Further, the deleterious effect of noise can also be appreciated from Fig. 3 (B), where it is noticeable that high-noise measurements require a larger number of iterations for the proposed algorithm to converge. Nonetheless, it is important to note that this apparent slower convergence rate is caused mainly by a larger number of  candidate basis functions (often with low magnitude coefficients) remaining in the library after the pre-processing step (Fig. 3 (B-C)). That is, the proposed pre-processing approach is less effective  when the noise standard deviation increases. These results emphasize that using static regression approaches (such as OLS which is employed in pre-processing) does not necessarily lead to discovery of the true dynamics, even if the recovered models exhibit high goodness-of-fit. Despite the increasing size of the candidate basis function library after pre-processing, Fig. 3 (C) shows that the correct number of terms was identified via \ours{} for each dynamical system considered for increasing noise values. Furthermore, convergence to the governing equations was achieved in just a few number of thresholding steps (6-8 iterations) as seen in Fig. 3 (D).

Of particular interest are stiff systems (arising commonly in  applications relating to learning chemical reaction kinetics and mechanisms), whose dynamics have slow and fast components resulting on the states evolving on different time scales. Learning stiff differential equations is challenging due to the high computational cost involved in solving said systems (due mainly to the smaller integration step size requirements), as well as due to ill-conditioning of the relevant gradient matrices \cite{kim2021stiff}. For the systems under consideration shown in Fig. 3 (A), the van der Pol oscillator and the Brusselator  exhibit stiff dynamic behavior for the selected model parameters, which can be handled by employing advances discretization strategies to prevent the aforementioned numerical stability issues. \ours{} uses orthogonal collocation on finite elements \cite{biegler2010nonlinear}, which is one of the highest order methods (i.e., having the lowest approximation error) and for which relatively large time step sizes are allowed even for stiff equations, also improving on the computational effort required for this class of systems. As can be seen from the results in Fig. 3 and Fig. 4, \ours{} shows good performance for discovering the system dynamics regardless of the underlying models stiffness and under increasing measurement noise for van der Pol and Brusselator dynamics.

Validation results for the discovered fundamental equations are shown in Fig. 4 (A-B). The first quantitative step for validating the identified equations is computing relevant performance metrics such as the ones commonly used for regression tasks (e.g. mean squared error, mean absolute error, coefficient of determination, etc. \cite{pedregosa2011scikit}). In the case of Fig. 4 (A), the true system dynamics were used to generate the reference trajectories to be used for validation. In practice, nonetheless, the original/smoothed training data  or a withheld testing set can be used to assess the  predictive power of the discovered governing equations. While the mean squared error between the trajectories from the discovered equations and that of the true dynamics increases with data measurement noise (Fig. 4 (A)), the MSE is generally low relative to the magnitude of the measured states for each dynamical system. Qualitatively, the discovered trajectories accurately reproduce the dynamic behavior of the true system as seen in the phase plane plots in Fig. 4 (B).  The exception to these results is the Lorenz oscillator, which due to its inherent chaotic nature is very sensitive to small perturbations to the model coefficients and initial conditions \cite{brunton2016discovering}. However, Fig. 4 (B) shows that  \ours{} successfully captures the attractor dynamics even if the trajectories simulated from the discovered equations do not overlap perfectly with the measurement data. Note that for all instances of the Lorenz oscillator \ours{} revealed the correct terms in the dynamics and resulted in coefficients estimates with small average percent error of $2.23\%$ and standard deviation of $0.14$ for the highest noise setting considered.

Differently from less transparent and interpretable ML frameworks,  \ours{} produces dynamic equations in which each functional term can be directly attributed to some underlying physical phenomena. In the case of the Lotka-Volterra equations, for example, the recovered terms in the dynamics can be given the following interpretations:
\begin{itemize}
\item For the discovered dynamics of prey population  $\dot{x}_1 = x_1 - 0.01x_1x_2$: The first term ($x_1$) reflects the fact that members of the prey species reproduce at an exponential rate in the absence of predators. The second term ($-0.01x_1x_2$) can be explained by the fact that the probability that members of the prey and predator species meet is proportional to the product of their populations, and so is the predation rate which drives a reduction in  the number of prey individuals.
\item For the discovered dynamics of predator population $\dot{x}_2 = -x_2 + 0.02x_1x_2$: The first term ($-x_2$) corresponds to the fact that  predators die or leave the ecosystem (emigrate) in the absence of prey, which leads to an exponential decay in their numbers. The second term ($0.02x_1x_2$), similar to the predation rate, can be explained by the fact that the predator population growth depends on the availability of prey, but this growth rate need not be the same rate at which the predators consume prey.
\end{itemize} This type of analysis is a critical step in validating the discovered models, thus any available domain expertise should be leveraged when available to ensure that \ours{} equations are explainable based on the observed behavior of the physical system.

Data-driven discovery of governing laws is a promising avenue for advancing our understanding of and elucidating new phenomena across a wide range of disciplines. This new fundamental knowledge can in turn be used to drive the development of new technologies to solve pressing human-centered challenges. In this paper, we introduced and validated \ours, a novel moving horizon, nonlinear dynamic optimization framework for learning governing equations from noise-contaminated state measurements over time.  \ours{} leverages a discretized model of the system dynamics to estimate the basis function coefficients, as opposed to prior works that use sparse regression techniques to predict the approximate values of the state derivatives as a linear combination of the basis functions evaluated at the the measured data. We demonstrated \ours's main advantages using a variety of dynamical systems including highly stiff nonlinear differential equations, and showed that \ours{} is highly robust to noise-contaminated data.

\section*{Materials and Methods}

\subsection*{\ours{} pre-processing}

\subsubsection*{Data smoothing}
Data are assumed to be contaminated with noise  (which we assume to be Gaussian and with zero mean but with unknown standard deviation). Data smoothing techniques are employed to reveal patterns otherwise hidden by noise, and to aid model training. To this end, the Savitzky-Golay filter (SVGF) \cite{savitzky1964smoothing} is employed, which is based on local least-squares regression polynomial approximations applied to the data on moving windows of a given size. Since data sets  with greater amounts of noise require more smoothing (accomplished by e.g. longer filter windows and lower order polynomials),  an iterative smoothing strategy is proposed to determine the appropriate window size for each of the measured state variables. The proposed smoothing scheme is detailed in SM Algorithm 1. The algorithm consists of increasing the Savitzky-Golay filter window size iteratively until there is no significant reduction in the standard deviation of the differenced initial measurement time series (in this sense, e.g. the first difference of a time series $\textbf{y}$ denoted $\textbf{y}'$  has entries given by  $\textbf{y}'(t) := \textbf{y}(t)-\textbf{y}(t-1)$). It is expected that data sets with noise of higher magnitude will take a greater number of iterations, and thus larger window sizes. An illustration of the results produced by the smoothing algorithm for data sets with different amounts of measurement noise collected from  the Lotka-Volterra system are shown in the SM Fig. S.1.

\subsubsection*{Statistical analysis}

Statistical tests are applied to the smoothed data to  infer which   features (i.e., basis functions) are most important for predicting the systems dynamics. Similar to \cite{brunton2016discovering}, we arrange the smoothed state measurements in the form of a data matrix $\tilde{\textbf{X}}$:
\begin{equation}
\label{eq:data_matrix}
\tilde{\textbf{X}}  =\begin{bmatrix}
\tilde{\textbf{x}}^{T}(t_1) \\ \vdots \\ \hat{\textbf{x}}^{T}(t_m)
\end{bmatrix}  = \begin{bmatrix}
\tilde{x}_1 (t_1) & \dots & \hat{x}_n (t_1)  \\
\vdots & \ddots & \vdots \\
\tilde{x}_1 (t_m) & \dots & \hat{x}_n (t_m) \\
\end{bmatrix}
\end{equation}
where $t_1$, $t_2$, $\dots$, $t_m$ are the time intervals at which the measurements were collected. The dictionary of candidate basis functions is then evaluated at every point in the data matrix. This can result in a structure such as the one below: 
\begin{equation}
\label{eq:basis_functions}
\Theta (\tilde{\textbf{X}} ) = \begin{bmatrix}
\vertbar  & \vertbar & \vertbar &  & \vertbar & &  \vertbar & & \vertbar \\
\textbf{1} & \tilde{\textbf{X}} & \tilde{\textbf{X}}^{P_2} & \cdots & 1/\tilde{\textbf{X}} & \cdots & e^{\tilde{\textbf{X}}} & \cdots & \sin (\tilde{\textbf{X}})   \\
\vertbar  & \vertbar & \vertbar &  & \vertbar & &  \vertbar & & \vertbar
\end{bmatrix}
\end{equation}
where $\tilde{\textbf{X}}^{P_2}$ denotes second order polynomials which may include interaction terms, e.g.:
\begin{equation}
\label{eq:second_order_poly}
\tilde{\textbf{X}}^{P_2}  = \begin{bmatrix}  \tilde{\textbf{x}}^{P_2}(t_1) \\ \vdots \\ \tilde{\textbf{x}}^{P_2}(t_m) \end{bmatrix}  = \begin{bmatrix}
\tilde{x}_1^2(t_1) & \tilde{x}_1(t_1) \tilde{x}_2(t_1)  & \cdots & \tilde{x}_2^2(t_1) & \cdots & \tilde{x}_n^2(t_1)\\
\tilde{x}_1^2(t_2) & \tilde{x}_1(t_2) \tilde{x}_2(t_2)  & \cdots & \tilde{x}_2^2(t_2) & \cdots & \tilde{x}_n^2(t_2) \\
\vdots & \vdots & \ddots & \vdots & \ddots & \vdots \\
\tilde{x}_1^2(t_m) & \tilde{x}_1(t_m) \tilde{x}_2(t_m)  & \cdots & \tilde{x}_2^2(t_m) & \cdots & \tilde{x}_n^2(t_m)
\end{bmatrix}
\end{equation}
Evidently, the form of this structure will depend on the choice of basis functions. The derivative of each state variable is approximated from the data by using central differences as follows:
\begin{equation}
\label{eq:central_difference}
\dot{\tilde{\textbf{x}}}(t_k) = \frac{\tilde{\textbf{x}}({t_{k+1}})-\tilde{\textbf{x}}({t_{k-1}})}{t_{k+1}-t_{k-1}}, \;\; k = 2,\dots,m-1
\end{equation}
and similar to \eqref{eq:basis_functions} a data derivative matrix $\dot{\tilde{\textbf{X}}}$ is formed. Other derivative approximation strategies (e.g., \cite{chartrand2011numerical}) can be used when there is significant noise in the data.

The first statistical test implemented is the Granger causality test \cite{granger1969investigating}, which establishes that: if a signal $\textbf{y}$ ``Granger-causes'' a signal $\textbf{z}$, then the past values of $\textbf{y}$ should contain information that helps predict $\textbf{z}$ above and beyond the information contained in the past values of $\textbf{z}$ alone. The \verb"statsmodels"
 \cite{seabold2010statsmodels} implementation of the Granger Causality test was used in Python. The test was designed to assess whether each of the basis functions evaluated on $\tilde{\textbf{X}}$ provides meaningful information in predicting future values of ${\tilde{\textbf{X}}}$ (i.e., predicting the evolution of future system states over time). In the typical formulation of the Granger causality test, the null hypothesis is $\Theta_i(\tilde{\textbf{X}}_j)$ does not Granger-cause ${\tilde{\textbf{X}}}_j$ for a given basis function $i\in\{1,\dots,n_\theta\}$ and a given state $j\in\{1,\dots,n_x\}$, which is evaluated by fitting autoregressive models of the form:
\begin{equation}
\label{eq:AR_models}
\begin{aligned}
 \tilde{\textbf{X}}_j(t) & = a_0 + a_1\tilde{\textbf{X}}_j(t-1) +\dots+a_m {\tilde{\textbf{X}}}_j(t-m)+e(t) \\
{\tilde{\textbf{X}}}_j(t) &= \hat{a}_0 + \hat{a}_1 {\tilde{\textbf{X}}}_j(t-1) + \dots+\hat{a}_m {\tilde{\textbf{X}}}_j (t-m)+b_1\Theta_i(\tilde{\textbf{X}}_j(t-1))+ \dots\\
&+b_q \Theta_i(\tilde{\textbf{X}}_j(t-q))+\hat{e}(t)
\end{aligned}
\end{equation}
where in this case we are only interested in a single lagged value of the evaluated basis functions (i.e.,$\Theta_i(\tilde{\textbf{X}}_j(t-1))$). The significance of using $\Theta_i(\tilde{\textbf{X}}_j)$ to predict ${\tilde{\textbf{X}}}_j$ is determined by examining the variance of the residuals $e(t)$ and $\hat{e}(t)$ by performing statistical tests based on F and chi-squared distributions. If the average p-value across all tests is less than the given significance level, then basis function $i$ is kept in the dictionary (i.e., the null hypothesis is rejected) for the state $j$. Otherwise, the basis function is removed from the dictionary for subsequent steps. The Granger causality test is repeated for every state variable and each basis function in our library. It should be noted that stationarity of the time series might need to be enforced (e.g., by differencing the signals and checking for stationarity via the Dickey-Fuller test \cite{dickey1979distribution}).

The second, and likely the most important, step in the proposed pre-processing statistical analysis involves an ordinary least-squares (OLS) regression problem to obtain a preliminary estimate of the coefficients $\Xi$. This step not only provides a solution to initialize DNLP embedded in the moving horizon discovery process, but also can also be used to derive lower and upper confidence bounds on $\Xi$, that are important in improving the convergence of the corresponding nonlinear programming problem. The linear model to be determined by performing OLS is given by:
\begin{equation}
\label{eq:OLS_preprocessing}
\dot{\tilde{\textbf{X}}} = \Theta(\tilde{\textbf{X}})\Xi^{OLS}
\end{equation}
where the coefficients $\Xi^{OLS} \in \mathbb{R}^{n_\theta \times n_x}$ are to be estimated. Note that a separate regression problem is solved for each of the states $1,\dots, n_x$ to compute each of the columns $\pmb{\xi}_i^{OLS}$ of $\Xi^{OLS}$. The OLS problems were implemented in the \verb"statsmodels" package \cite{seabold2010statsmodels} in Python, which leverages linear algebra tools to efficiently solve  the normal equations to estimate the coefficient vector for state $i$:
\begin{equation}
\label{eq:normal_eq}
\pmb{\xi}_i^{OLS} = (\Theta(\tilde{\textbf{X}})^T\Theta(\tilde{\textbf{X}}))^{-1}(\Theta(\tilde{\textbf{X}})\dot{\tilde{\textbf{X}}}_i
\end{equation}
where $\dot{\tilde{\textbf{X}}}_i$ represents the  approximation of the derivative of state $i$ estimated from the data.

As a preliminary approach to selecting the most informative basis functions, the results obtained by OLS can be used to intuit  the most important predictor variables within $\Theta(\tilde{\textbf{X}})$. In brief, classical approaches for performing this variable selection procedure involve computing the F-statistic and examining the corresponding p-values \cite{james2013introduction}, which are automatically computed by solving \eqref{eq:OLS_preprocessing} using the \verb"statsmodels" OLS implementation \cite{seabold2010statsmodels}. In this sense, under the null hypothesis that a coefficient for a basis function is zero, predictors having p-values greater than a specified significance level are eliminated. For the coefficients with p-values small enough that the null hypothesis cannot be confidently rejected, we also use the OLS results to derive confidence intervals (for a pre-specified confidence level) in order to obtain lower ($\Xi^L$) and upper ($\Xi^U$)  bounds on the coefficient values. For a more extensive discussion on the statistical properties of OLS, the reader is referred to established texts on  machine learning \cite{bishop2006pattern, james2013introduction}, as well as the \verb"statsmodels" documentation \cite{seabold2010statsmodels}.

While these steps can be a useful preliminary approach to eliminating basis functions from the dictionary in a specific application and for a given data set, they inherently rely  on estimating the derivative from noisy state measurements (which introduces additional error to the already noise-contaminated data). To this end, we suggest that a high significance level should be used (in both Granger causality tests and OLS regression) for eliminating basis functions from the dictionary; this is a conservative approach that prevents eliminating the basis functions that do belong in the true system dynamics. Further, since high-noise measurements may significantly affect the accuracy of coefficient estimates in OLS, it is recommended that large confidence intervals (e.g., 99.99\%) be  used to estimate $\Xi^L$ and $\Xi^U$.

\subsection*{\ours{} dynamic nonlinear optimization}

\subsubsection*{Discretization of the dynamic equations}

In this work, discretization is performed with respect to the time domain (i.e., by defining a finite set of discrete points in time where the dynamic equations are evaluated), to convert the model in continuous time in \eqref{eq:dyn_basis} to a discrete time expression of the form of $\textbf{x}({k+1}) = \textbf{g}(\Theta(\textbf{x}(k)),\Xi)$ as used  in optimization problem \eqref{eq:dyn_opt}. Simultaneous strategies  \cite{biegler2010nonlinear} are employed, whereby the discretized equations are incorporated as nonlinear algebraic constraints in problem \eqref{eq:dyn_opt}. One of the simplest and most widely used class of methods are finite difference transformations, such as explicit and implicit Euler schemes. While generally more challenging to implement, collocation methods provide  substantially more accurate approximations of the dynamical system, have good numerical stability allowing relatively large time steps to be considered, and are thus advantageous for stiff dynamical systems \cite{biegler2010nonlinear}. Broadly speaking, collocation on finite elements entails partitioning the time domain into $M-1$ finite elements, over which polynomials of order $K+1$ are used to approximate the differential variable $\textbf{x}(t)$ (each polynomial for each finite element is defined using $K$ collocation points,  which act as an additional discretization within each finite element). Additional constraints are introduced to enforce continuity across the finite element boundaries for each differential variable: \cite{biegler2010nonlinear}.
\begin{equation}
\begin{aligned}
\label{eq:collocation1}
& \left. \frac{d\textbf{x}(t)}{dt} \right \vert_{t_{ij}} =  \frac{1}{h_i}\sum_{j=0}^{K}\textbf{x}_{ij}\frac{d\ell_j(\tau_k)}{d\tau}, \;\; k \in \{1,\dots,K\}, \; i\in \{ 1,\dots,M-1\} \\
& \textbf{x}_{i+1,0} = \sum^{K}_{j=0} \ell_j(1)\textbf{x}_{i,j}, \;\; i\in \{1,\dots,M-1\} \\
\end{aligned}
\end{equation}
where the state variable $\textbf{x}(t)$ is interpolated using Lagrange polynomials as follows:
\begin{equation}
\begin{aligned}
\label{eq:collocation}
& t_{ij} = t_{i-1}+\tau_jh_i \\
& \textbf{x}(t) = \sum^{K}_{j=0}\ell_j(\tau)x_{ij}, \;\; t\in [t_{i-1}, t_i], \; \tau \in [0,1] \\
& \ell_j(\tau) = \prod^{K}_{k=0,\neq j} \frac{\tau-\tau_k}{\tau_j-\tau_k}
\end{aligned}
\end{equation}
A comprehensive discussion on discretization strategies can be found elsewhere \cite{biegler2010nonlinear}. In this work, we leverage the \verb"pyomo.DAE" \cite{nicholson2018pyomo} modeling extension that enables automatic simultaneous discretization of ODEs, and leverages Gauss-Legendre and Gauss-Radau collocation schemes to determine the interpolating points in \eqref{eq:collocation}. . We note that the collocation equations in \eqref{eq:collocation1} and \eqref{eq:collocation} used to discretize the dynamics in the form of \eqref{eq:dyn_basis} are represented in compact form as $\textbf{x}({k+1}) = \textbf{g}(\Theta(\textbf{x}(k)),\Xi)$ in the DNLP in \eqref{eq:dyn_opt}.

It should be noted that  the optimal choice of collocation points may not align with the sample times  $t_1, \dots, t_m$ at which the data were originally collected, thus spline interpolation is required to approximate the data at the relevant time instants. We employ the \verb"SciPy" package \cite{virtanen2020scipy} in Python using a cubic spline to estimate the state values at the collocation points.

\subsection*{\ours{} moving horizon optimization}

While the formulation introduced in \eqref{eq:dyn_opt} has several advantageous properties relative to prior regression-based framework, a key challenge is addressing the computational burden of solving this problem as the size of the data set and the dimension of the dictionary of basis functions increase. Longer time horizons for the data set require  a larger number of finite elements and collocation points. This in turn increases the number of variables and constraints in \eqref{eq:dyn_opt}. To alleviate the computational burden of solving the optimization problem in \eqref{eq:dyn_opt}, we employ ideas stemming from control and estimation theory. In particular, we draw inspiration from moving horizon estimation (MHE) \cite{rao2003constrained}, which involves solving a sequence of state estimation problems (typically in the form of a DNLP) online and discarding old measurements for which state estimates have already been computed (instead of estimating states, \ours{} estimates the coefficients $\Xi$). For the purposes of \ours, these operations need not be performed online; the idea is that the original data set can be segmented into several sequential subsets of smaller size over which different instances of \eqref{eq:dyn_opt} are solved.

The key elements of the moving horizon strategy are illustrated in Fig. S.2 in the SM for data corresponding to the Lotka-Volterra system. Note that the number of optimization problems to be solved depends directly on the choice of optimization horizon $H$. Nevertheless, the problems in the sequence are likely significantly more computationally tractable than solving  \eqref{eq:dyn_opt} for the entire (large-scale) data set. It is worth mentioning that, to date, no analytical frameworks exist for determining the optimal choice of horizon $H$. The empirical consensus is that longer  horizons yield better results (i.e., convergence of the estimates to the true values of the parameters), which intuitively comes at a computational cost \cite{rao2002constrained}. For periodic systems, such as the Lotka-Volterra system shown in SM Fig. S.2, an intuitive choice for $H$ can be an integer multiple of the period corresponding to the fundamental oscillation frequency, which can be estimated from the data.

A detailed outline of the  moving horizon algorithm is presented next in SM  Algorithm S.2. In brief, the algorithm consists of analyzing the set of smoothed state measurements in a moving horizon fashion as illustrated in SM Fig S.2., solving \eqref{eq:dyn_opt} and performing parameter thresholding every $\omega$ iterations (the thresholding process is described in Algorithm S.3).  The moving horizon algorithm terminates when  the training data are exhausted (recall that the original data set consists of a total of $m$ measurements collected at times $t_1,\dots,t_m$) or convergence is established, that is when the number of basis functions remaining in the library, denoted as $|\Theta|$, does not change after a number $\Omega$  of thresholding steps. Upon convergence, the outputs of the algorithm are set of discovered basis function and the corresponding mean values of the coefficients obtained for the last $\omega\times\Omega$ iterations, for which $|\Theta|$ did not change. If the algorithm fails to converge and the data set is exhausted, the following additional steps can be attempted: collecting a larger data set, extending the dictionary of basis functions -- possibly leveraging domain knowledge, and/or deriving tighter coefficient bounds.

\subsubsection*{Thresholding algorithm}

A key component of \ours{} is systematically eliminating non-basic functions from the initial dictionary. The proposed thresholding approach is outlined in Algorithm S.3, and is embedded within the moving horizon  scheme introduced previously in Algorithm S.2. The proposed framework leverages the fact that the values of non-basic coefficients are expected to be low in magnitude but to have significant variability when estimated using different subsets of the data. Such effects can be quantified statistically by computing the coefficient of variation ($CV$), defined as the ratio between the standard deviation and the mean, evaluated for a series of coefficient estimates obtained from successive portions of the data. In our framework, the coefficient of variation is computed for the coefficients of all basis functions $\Theta_\theta \; \forall \theta \in \{1,\dots,|\Theta |\}$, and for all state variables $j\in\{1,\dots,n_x\}$. If the coefficient of variation $CV_{\theta,j}$ is greater than the specified variability threshold $\gamma$, then the basis function is pruned and not considered in future iterations of the moving horizon scheme. Otherwise, the basis function $\Theta_\theta$ remains in the dictionary.

It is important to note that some basis functions are significantly more susceptible to thresholding that others. For example, basis functions such as  $\{\pmb{1}, \textbf{x}\} \in \Theta$ are particularly prone to contributing to overfitting the initial measurement noise in the data, as well as the differentiation error introduced when the dynamics are discretized (particularly in high noise environments and when the initial function library is larger). To prevent spurious thresholding for this type of basis functions whose associated coefficients are likely so see greater variability across different data subsets, an alternative is to keep them in the basis for the first few thresholding steps regardless of their associated coefficient's observed $CV$. After these initial iterations and when (potentially) some of the other non-basic basis functions have been pruned, if basis functions like e.g. $\{\pmb{1}, \pmb{x}\}$ are in fact basic they are expected to experience less variability in their respective coefficients and remain in the basis when the algorithm converges.

\bibliographystyle{unsrt}
\bibliography{DySMHO_paper_v6_bib}

\section*{Acknowledgments:}

\subsection*{Funding:}
Support from the National Science Foundation, USA through the CAREER Award 1454433 (recipient: MB) is acknowledged with gratitude.

\subsection*{Competing interests}
The authors declare that they have no competing interests.

\subsection*{Data and materials availability}
All data and code used in this analysis can be found at: \url{https://github.com/Baldea-Group/DySMHO}

\clearpage 

\section*{Supplementary Materials}

\renewcommand{\thesection}{S\arabic{section}}
\renewcommand{\thesubsection}{S\arabic{subsection}}
\renewcommand{\theequation}{S.\arabic{equation}}
\renewcommand{\thefigure}{S.\arabic{figure}}
\renewcommand{\thetable}{S.\arabic{table}}
\renewcommand{\thealgocf}{S.\arabic{algocf}}

\subsection[S]{\ours{} pre-processing}

\vspace{0.5cm}
\begin{algorithm}[H]
\SetAlgoLined
\KwResult{smoothed state measurement data $\tilde{\textbf{x}}$}
 \textbf{Inputs:} raw state measurement data $\hat{\textbf{x}}$, initial window size $WS$, step size $\Delta$, polynomial order $\phi$, smoothing threshold $\alpha$\;
Standardize measured data using min-max scaling : $\hat{\textbf{x}}'(t)=\frac{\hat{\textbf{x}}(t)-\min_t \hat{\textbf{x}}(t) }{\max_t \hat{\textbf{x}}(t) - \min_t \hat{\textbf{x}}(t)}$ \;
Calculate differenced state measurement as proxy for noise: $\hat{\textbf{x}}''(t) = \hat{\textbf{x}}'(t) - \hat{\textbf{x}}'({t-1}) $\;
 Calculate standard deviation of noise: $\sigma_{i} \leftarrow \sigma(\hat{\textbf{x}}'')$\;
 $WS_i \leftarrow 0, WS_j \leftarrow WS$ \;
 \While{\textit{True}}{
Perform smoothing: $\tilde{\textbf{x}'} = \text{SVGF}(\hat{\textbf{x}}',WS_j,\phi)$\;
 $\tilde{\textbf{x}}''(t) = \tilde{\textbf{x}}'(t) - \tilde{\textbf{x}}'({t-1}), \sigma_j \leftarrow \sigma(\tilde{\textbf{x}}'')$ \;
  \eIf{$|\sigma_i -\sigma_j|/|\sigma_i| <\alpha$}{
 \Return $\tilde{\textbf{x}} \leftarrow \text{SVGF}(\hat{\textbf{x}},WS_i,\phi)$
   }{
   Increase filter window size: $WS_i \leftarrow WS_j, WS_j \leftarrow WS_j+\Delta$
  }
 }
 \caption{Iterative smoothing via Savitzky-Golay filtering}
\end{algorithm}

\begin{figure}[h!]
\centering
\includegraphics[scale = 0.5]{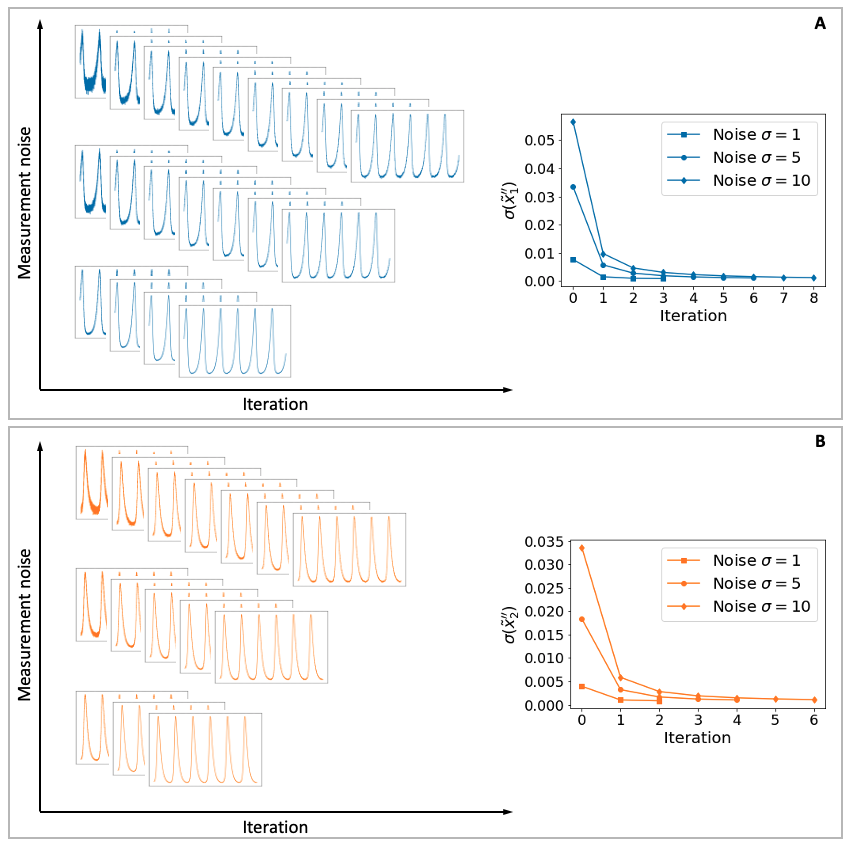}
\caption{Illustration of smoothing algorithm for Lotka-Volterra predator-prey model. (\textbf{A})  results obtained for state $x_1$. (\textbf{B})  results for state $x_2$. Figures on the left show the smoothed system trajectories as a function of initial measurement noise at every iteration of the smoothing algorithm, figures on the right show the standard deviation of the estimated noise  at every iteration. (Results obtained using $\alpha=0.1$, $WS_i=10$, $\Delta=10$, and $\phi = 2$)}
\end{figure}

\clearpage

\subsection[S]{\ours{} moving horizon optimization}

For each data subset the following nonlinear optimization problem is solved:
\begin{equation}
\label{eq:dyn_opt_disc}
\begin{aligned}
& \min_{\Xi} && \frac{1}{2M} \sum_{k=1}^{M}||\textbf{x}(k)-\tilde{\textbf{x}}(k)||_2^2 + \lambda \ell(\Xi)\\
& \; \text{s.t.} && \textbf{x}({k+1}) = \textbf{g}(\Theta(\textbf{x}(k)),\Xi)  \;\;\; \forall k\in \{1,\dots,M\}\\
& \; &&\textbf{x}(1) = \tilde{\textbf{x}}(t_1) \\
& \; && \Xi \in \{ \Xi^{L}, \Xi^{U} \}
\end{aligned}
\end{equation}

\begin{figure}[h!]
\centering
\includegraphics[scale = 0.25]{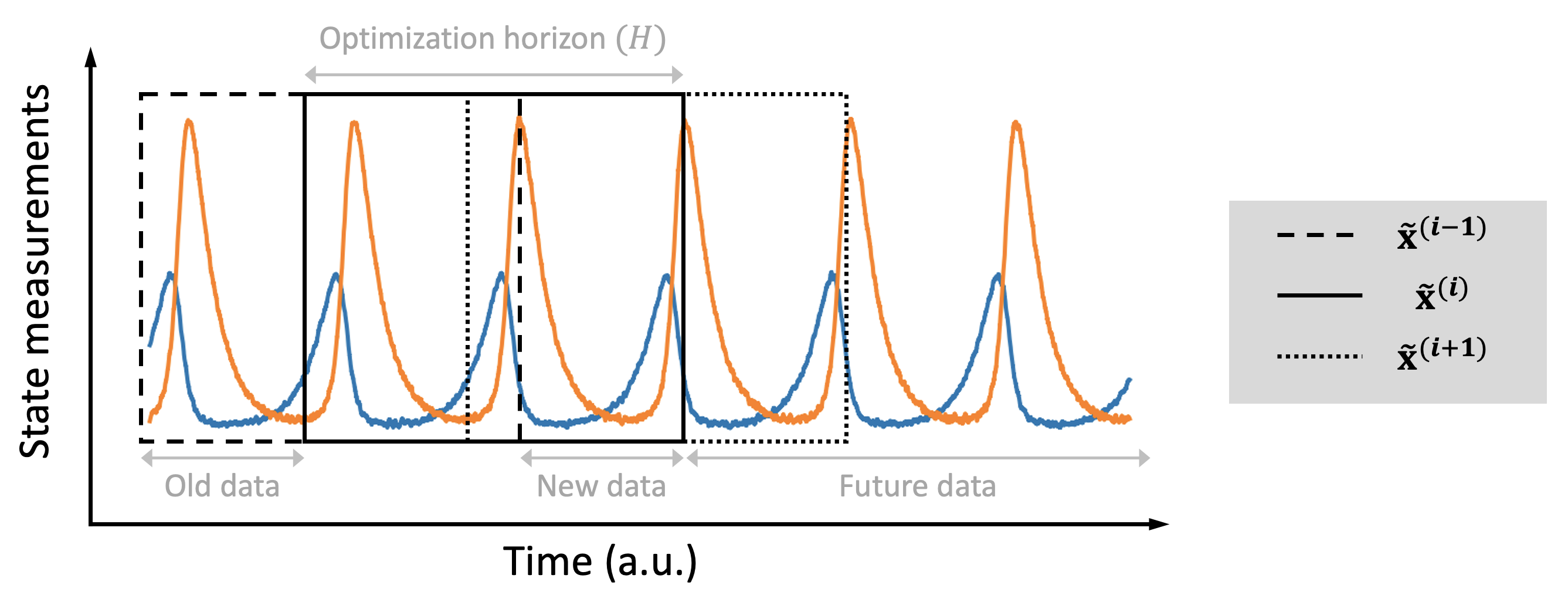}
\caption{Illustration of subsets of data of length $H$ used for the proposed moving horizon algorithm using simulated and smoothed state measurements corresponding to the Lotka-Volterra predator-prey model. Data frames used for the past iteration, current iteration, and next iteration are denoted as $\tilde{\textbf{x}}^{(i-1)}$, $\tilde{\textbf{x}}^{(i)}$, and $\tilde{\textbf{x}}^{(i+1)}$, respectively.
}
\end{figure}

\begin{algorithm}[H]
\SetAlgoLined
\KwResult{underlying governing dynamics given by $\dot{\textbf{x}}=\Xi^T\Theta(\textbf{x}^T)^T$ }
 \textbf{Inputs:} smoothed state measurement data $\tilde{\textbf{x}}$, initial library of basis functions $\Theta(\cdot)$, optimization horizon $H$, data step size $\Delta \mathcal{D}$, thresholding frequency $\omega$, tresholding steps for convergence  $\Omega$, $\Xi^L, \Xi^U$\;
 Set counters $i\leftarrow 1$, $j \leftarrow 0$ \;
 \While{$H+(i-1)\Delta \mathcal{D} \le m$}{
Get data for current iteration: $\tilde{\textbf{x}}^{(i)}\leftarrow \tilde{\textbf{x}}(i\Delta \mathcal{D}:H+i\Delta \mathcal{D})$\;
Obtain coefficient estimates $\Xi^{(i)}$ by solving \eqref{eq:dyn_opt_disc} using $\tilde{\textbf{x}}^{(i)}, \Theta(\cdot), \Xi^L, \Xi^U$\;
  \If{$\sim i \% \omega $}{
  Perform thresholding step described in Algorithm S.3 to compute $\bar{\Theta}$  \;
  \eIf{$|\Theta| = |\bar{\Theta}|$}{
  $j\leftarrow j+ 1$ \;
  \If{$j =\Omega$}{
   \Return $\Theta$, $\Xi \leftarrow \frac{1}{\omega\times \Omega}\sum_{k=i-\Omega}^{k=\Omega} \Xi^{(k)}$ }
  }{
    Recompute $\Xi^L, \Xi^U$ via OLS for the reduced basis function library $\bar{\Theta}$ \;
  Set $\Theta \leftarrow \bar{\Theta}$ \;
  }
  }
   Set $i \leftarrow i+1$,
 }
 \caption{Moving horizon optimization algorithm}
\end{algorithm}

\clearpage

\subsection[S]{\ours{} thresholding algorithm}

\begin{algorithm}[H]
 \label{alg:thresholing}
\SetAlgoLined
\KwResult{updated library of basis functions $\bar{\Theta}(\cdot)$}
 \textbf{Inputs:} current library of basis functions $\Theta$, coefficient estimates for iterations $\Xi^{(i)},\Xi^{(i+1)}, \dots, \Xi^{(i+\omega-1)},\Xi^{(i+\omega)}$, variability tolerance $\gamma$\;
 \For{$\theta \in \{1,\dots,|\Theta|\}$}{
 \For{$j \in \{1,\dots,n_x\}$}{
 Compute average: $\mu_{\theta,j} \leftarrow \frac{1}{\omega}\sum_{k=i}^{k=i+\omega}\Xi_{\theta,j} ^{(k)}$\;
  Compute standard deviation: $\sigma_{\theta,j}  \leftarrow \frac{1}{\omega}\sum_{k=i}^{k=i+\omega} (\Xi_{\theta,j} ^{(k)}-\mu_{\theta,j} )^2$ \;
  Compute coefficient of variation: $CV_{\theta,j} \leftarrow \sigma_j / \mu_j$\;
  \If{$CV_{\theta,j}  < \gamma $ }{
  Append $\Theta_\theta$ to $\bar{\Theta}$
  }
 }
 }
 \Return $\bar{\Theta}$
 \caption{Thresholding algorithm to prune the dictionary of candidate basis functions}
\end{algorithm}

\clearpage

\subsection[S]{\ours{} configurations for numerical experiments}

All computations were performed on a PC running Windows 7 64-bit, with a 3.6 GHz Intel Core i7-7700 processor and 32 GB RAM. All DNLPs where formulated in Python 3.8.3 using \verb"Pyomo" \cite{hart2017pyomo}, and solved using CONOPT \cite{drud1994conopt} as the nonlinear solver using all default settings.  The objective function minimized was the mean squared difference between the model and the state measurements, without regularization (i.e., $\lambda = 0$). The basis functions for the Lotka-Volterra, van der Pol, and Brusselator examples were:
\begin{equation}
\begin{aligned}
\label{eq:basis_fun_2D}
& \Theta_{x_1} = \{1,x_1,x_2,x_1x_2,x_1^2,x_2^2,x_1^2x_2,x_1x_2^2,x_1^3,x_1^4,1/x_1,e^{x_1}, \sin{x_1}, \cos{x_1} \} \\
& \Theta_{x_2} = \{1,x_2,x_1,x_1x_2,x_2^2,x_1^2,x_2^2x_1,x_2x_1^2,x_2^3,x_2^4,1/x_2,e^{x_2}, \sin{x_2}, \cos{x_2} \}
\end{aligned}
\end{equation}
and for the Lorenz example were:
\begin{equation}
\begin{aligned}
\label{eq:basis_fun_3D}
& \Theta_{x_1} = && \{1,x_1,x_2,x_3,x_1x_2,x_1x_3,x_2x_3, x_1^2, x_2^2, x_3^2, x_1^2x_2, x_1x_2^2, x_1^2x_3, x_1x_3^2, x_2^2x_3,  \\&  && x_2x_3^2, x_1^3, x_1^4, 1/x_1, e^{x_1}, \sin{x_1}, \cos{x_1} \}  \\
& \Theta_{x_2} = && \{1,x_2,x_1,x_3,x_2x_1, x_2x_3, x_1x_3, x_2^2, x_1^2, x_3^2, x_2^2x_1, x_2x_1^2, x_2^2x_3, x_2x_3^2, x_1^2x_3,   \\&  && x_1x_3^2, x_2^3, x_2^4, 1/x_2, e^{x_2}, \sin{x_3}, \cos{x_3} \}  \\
& \Theta_{x_3} = && \{1,x_3,x_1,x_2,x_3x_1, x_3x_2, x_1x_2, x_3^2, x_1^2, x_2^2, x_3^2x_1, x_3x_1^2, x_3^2x_2, x_3x_2^2, x_1^2x_2,   \\&  && x_1x_2^2, x_3^3, x_3^4, 1/x_3, e^{x_3}, \sin{x_3}, \cos{x_3} \}  \\
\end{aligned}
\end{equation}
Tables S.1-S.4 outline all of the \ours{} parameters and configurations used for the numerical experiments.

\begin{table}[htb]
\caption{\ours{} configuration for Lotka-Volterra predator-prey model numerical experiments}
\vspace{-0.2in}
\label{tb:LV_dysmo}
\begin{center}
\begin{tabular}{l l l}
\hline
\textbf{System}  &  \textbf{\ours{} configuration} &  \\
\hline
\hline
\textit{Lotka-Volterra}  & {Data simulation} & Initial conditions: (100,15)	\\
 & & Sampling frequency: 1/500  \\
 &  {Smoothing} 	& $WS=10$   \\
 & & $\gamma=10$   \\
 & & $\alpha=0.1$   \\
 & {Pre-processing} & Granger tests: $\chi^2$ and F-distributions \\
 &  & Granger p-value: 0.1 \\
 & & OLS p-value: 0.9 \\
 & & OLS \% confidence: $1\times 10^{-6}$ \\
 & {Discretization} & Scheme: Lagrange-Radau \\
 & & Finite elements: 50 \\
 & & Collocation points: 15 \\
 & & Data interpolation: Cubic spline  \\
 & {Moving horizon} & $H$: 6 (3,000 data samples) \\
 & & $\Delta \mathcal{D}$: 100 samples \\
 & & $\omega$: 10 \\
 & & $\Omega$: 40 \\
 & {Thresholding} & $\gamma$: 1 \\
\hline
\end{tabular}
\end{center}
\end{table}

\begin{table}[htb]
\caption{\ours{} configuration for van der Pol oscillator numerical experiments}
\vspace{-0.2in}
\label{tb:VDP_dysmo}
\begin{center}
\begin{tabular}{l l l}
\hline
\textbf{System}  &  \textbf{\ours{} configuration} &  \\
\hline
\hline
\textit{Van der Pol}  & {Data simulation} & Initial conditions: (1,-2)	\\
 & & Sampling frequency: 1/500  \\
 &  {Smoothing} 	& $WS=10$   \\
 & & $\gamma=10$   \\
 & & $\alpha=0.1$   \\
 & {Pre-processing} & Granger tests: $\chi^2$ and F-distributions \\
 &  & Granger p-value: 0.1 \\
 & & OLS p-value: 0.8 \\
 & & OLS \% confidence: $1\times 10^{-6}$ \\
 & {Discretization} & Scheme: Lagrange-Radau \\
 & & Finite elements: 80 \\
 & & Collocation points: 15 \\
 & & Data interpolation: Cubic spline  \\
 & {Moving horizon} & $H$: 20 (40,000 data samples) \\
 & & $\Delta \mathcal{D}$: 50 samples \\
 & & $\omega$: 10 \\
 & & $\Omega$: 40 \\
 & {Thresholding} & $\gamma$: 1 \\
\hline
\end{tabular}
\end{center}
\end{table}

\begin{table}[htb]
\caption{\ours{} configuration for Brusselator numerical experiments}
\vspace{-0.2in}
\label{tb:B_dysmo}
\begin{center}
\begin{tabular}{l l l}
\hline
\textbf{System}  &  \textbf{\ours{} configuration} &  \\
\hline
\hline
\textit{Brusselator}  & {Data simulation} & Initial conditions: (1,1)	\\
 & & Sampling frequency: 1/1000  \\
 &  {Smoothing} 	& $WS=10$   \\
 & & $\gamma=10$   \\
 & & $\alpha=0.1$   \\
 & {Pre-processing} & Granger tests: $\chi^2$ and F-distributions \\
 &  & Granger p-value: 0.1 \\
 & & OLS p-value: 0.8 \\
 & & OLS \% confidence: $1\times 10^{-6}$ \\
 & {Discretization} & Scheme: Lagrange-Radau \\
 & & Finite elements: 60 \\
 & & Collocation points: 15 \\
 & & Data interpolation: Cubic spline  \\
 & {Moving horizon} & $H$: 10 (10,000 data samples) \\
 & & $\Delta \mathcal{D}$: 100 samples \\
 & & $\omega$: 10 \\
 & & $\Omega$: 40 \\
 & {Thresholding} & $\gamma$: 1 \\
\hline
\end{tabular}
\end{center}
\end{table}

\begin{table}[htb]
\caption{\ours{} configuration for Lorenz chaotic oscillator numerical experiments}
\vspace{-0.2in}
\label{tb:L_dysmo}
\begin{center}
\begin{tabular}{l l l}
\hline
\textbf{System}  &  \textbf{\ours{} configuration} &  \\
\hline
\hline
\textit{Lorenz}  & {Data simulation} & Initial conditions: (-8,8,27)	\\
 & & Sampling frequency: 1/1000  \\
 &  {Smoothing} 	& $WS=10$   \\
 & & $\gamma=10$   \\
 & & $\alpha=0.1$   \\
 & {Pre-processing} & Granger tests: $\chi^2$ and F-distributions \\
 &  & Granger p-value: 0.1 \\
 & & OLS p-value: 0.7 \\
 & & OLS \% confidence: $1\times 10^{-6}$ \\
 & {Discretization} & Scheme: Lagrange-Radau \\
 & & Finite elements: 50 \\
 & & Collocation points: 15 \\
 & & Data interpolation: Cubic spline  \\
 & {Moving horizon} & $H$: 2 (2,000 data samples) \\
 & & $\Delta \mathcal{D}$: 100 samples \\
 & & $\omega$: 10 \\
 & & $\Omega$: 40 \\
 & {Thresholding} & $\gamma$: 1 \\
\hline
\end{tabular}
\end{center}
\end{table}

\end{document}